\documentclass[11pt,twoside]{article}
\usepackage{latexsym}
\usepackage{amssymb,amsbsy,amsmath,amsfonts,amssymb,amscd}
\usepackage{epsfig, graphicx}
\newcommand{\fer}[1]{(\ref{#1})}
\newcommand {\wt}[1] {{\widetilde #1}}
\newcommand {\al} {\alpha}
\newcommand {\e}  {\varepsilon}

\newcommand {\sg} {\sigma}
\newcommand {\f}   {\frac}
\newcommand {\p}   {\partial}
\newcommand{\dis}{\displaystyle}

\newcommand{\beq}{\begin{equation}}
\newcommand{\beqa}{\begin{eqnarray}}
\newcommand{\bea} {\begin{array}{ll}}
\newcommand{\beqan}{\begin{eqnarray*}}
\newcommand{\eeq}{\end{equation}}
\newcommand{\eeqa}{\end{eqnarray}}
\newcommand{\eeqan}{\end{eqnarray*}}
\newcommand{\eea} {\end{array}}
\newtheorem{theorem}{Theorem}[section]
\newtheorem{lemma}[theorem]{Lemma}

\newtheorem{corollary}[theorem]{Corollary}
\newcommand{\qed}{{ \hfill
                       {\unskip\kern 6pt\penalty 500
                       \raise -2pt\hbox{\vrule\vbox to 6pt{\hrule width 6pt
                       \vfill\hrule}\vrule} \par}   }}
\title{\Large \bf Asymmetric potentials and motor effect:  a large deviation approach}

\author{ Beno\^ \i t Perthame$^{1}$ and Panagiotis E. Souganidis$^{2,4}$}
\date{\today}

\begin{document}
\maketitle
\vspace*{1.0cm}
\pagenumbering{arabic}

\begin{abstract} We provide a mathematical analysis of appearance of the concentrations
(as Dirac masses) of the solution to a  Fokker-Planck system with asymmetric potentials. 
This problem has been proposed as a model to describe motor proteins moving along molecular 
filaments. The components of the system describe the densities of the different conformations of the proteins. 

Our results are based on the study of
a Hamilton-Jacobi equation arising, at the zero diffusion limit, after an exponential 
transformation change of the phase function that rises a Hamilton-Jacobi equation. We consider
different classes of conformation transitions coefficients (bounded, unbounded and locally 
vanishing).   

\end{abstract}
\pagestyle{plain} \vspace*{1.0cm}
\\
\noindent {\bf AMS Class. Numbers.} 35B25, 49L25, 92C05

\section{Introduction}
\label{Sec:intro}

\noindent A striking feature of living cells is their ability to generate motion, as, for 
instance in muscle contraction already investigated theoretically in the 50's (\cite{hux}). 
But even more elementary processes allow for intra-cellular material transport along various 
filaments that are part of the cytoskeleton. These are known as ``motor proteins''. For
example, myosins move along actin filaments and kinesins and dyneins move along micro-tubules. 
In the early 90's, it became possible to device a new generation of experiments 
{\em in vitro} where both the filaments and the motor proteins are sufficiently purified. This
lead to an improved biophysical understanding of the biomotor process  (see, for
instance, \cite{ah, jap, PEO, DEO}, and the tutorial book \cite{howard}) and gave rise to 
a large cellular biology literature. The experimental observations 
made possible to explain how chemical energy can
be transformed into mechanical energy and to come up with mathematical models for 
molecular motors. The underlying principles are elementary and represent in fact the common 
basis for all biomotors. On the one hand, the filament provides for an asymmetric potential
(this notion was introduced in the earliest theoretical descriptions by Huxley, \cite{hux}), 
sometimes referred to as the energy landscape.  On the other hand, the protein can reach 
several different conformations. This can
be ATP/ADP hydrolysis but five to six different states of the protein could be involved 
during muscular contraction.
\\

In this paper we consider the following model: Molecules can reach
$I$ configurations with density, for each $i=1,2...,I$, $n_i$. A bath
of such molecules is moving in an asymmetric potential seen differently by the $I$ 
configurations denoted, for $i=1,...,I$, by $\psi_i$.
Fuel consumption triggers a configuration change among the different states with rates 
$\nu_{ij} >0$, for $i,j=1,2...,I$. Diffusion, denoted below by $\sigma$, is taken
into account.
\\

Thes simple considerations lead to the following system of
elliptic equations for the densities $(n_i)_{1\leq i\leq I}$:
\beq
\left\{
\begin{array}{ll}
-\sg \f{\p^2}{\p x^2}  n_i - \f{\p}{\p x}( \nabla \psi_i\; n_i)
+\nu_{ii} n_i  =\dis{\sum_{j\neq i}}\; \nu_{ij} n_j  \qquad \text{
in } (0,1),
\\
\sg \f{\p}{\p x}  n_i (x)+ \nabla \psi_i (x) \; n_i (x)=0 \quad \text{ for $x=0$ or $1$.  }
\quad 
\end{array}
\right.
\label{eq:biomotor}
\eeq

The zero flux boundary conditions means that the total number of molecules, in each molecular state, is preserved by transport (but not by configuration exchange). 

Throughout the paper we assume that, for $i=1,...,I$
\beq n_i>0 \qquad \text{ in } [0,1].
\label{eq:positive}
\eeq

The zero flux boundary condition, motivated by the additional
modeling assumption that total density is conserved, leads to the
condition that, for all $i=1,...,I$,
\beq \label{as:nu}
\nu_{ii} = \sum_{j, \;j\neq i} \nu_{ji}.
\eeq

Several biomotor models, including the one described above, were
analyzed in \cite{chk,ckk,kk, hkm} through optimal transportation methods. In \cite{chk} it is proved that
there is a positive steady state solution that can, for instance, be  normalized by
\beq \label{eq:norm}
\int_0^1 \sum_{1\leq i \leq I} n_i(x) dx=1. \eeq

The simplest way to explain this fact is to observe that the adjoint system,
\beq
\left\{
\begin{array}{ll}
-\sg \f{\p^2}{\p x^2}  \phi_i + \nabla \psi_i\;  \f{\p}{\p x} \phi_i
+\nu_{ii} \phi_i =\dis \sum_{j\neq i} \nu_{ji} \phi_j   \qquad \text
{ in } (0,1),
\\
\\
\f{\p}{\p x}  \phi_i =0 \quad \text{ in } \quad \{0,1\}, 
\end{array}
\right. \label{eq:biomotoradj} \eeq
admits the trivial solution
$\phi_1=\phi_2=...=\phi_I=1$. This yields that $0$ is  the first
eigenvalue of the system and thus of its adjoint \fer{eq:biomotor}.
The Krein-Rutman theorem gives the $n_i$'s, but the solution is not explicitly known 
except for $I=1$, a situation  where the motor effect cannot be achieved. The stability 
of this problem is also related to the notion of relative entropy \cite{dkk,Pe,Pebook, mmp}.
\\

The typical results obtained about biomotors in
\cite{chk, hkm} are that, for small diffusion $\sg$, under some precise
asymmetry assumptions on the potentials, 
 the solutions tend to concentrate, as $\sg \to 0$, as Dirac masses at either $x=0$ or
$x=1$. In the sequel such a behavior will be called {\em motor effect}. 

Our results (i) provide an alternative proof of this motor effect, and (ii) allow for 
more general assumptions like, for instance, various scalings on the coefficients $\nu_{ij}$.
While \cite{chk, hkm} transform the system  \fer{eq:biomotor} into an ordinary differential 
equation and analyze directly its solution, here we use a direct PDE argument based on the 
phase functions $R_i=- \sg \ln n_i$ that satisfy (in the viscosity sense, \cite {BCD, barles,
 CIL, FS}) a Hamilton-Jacobi solution. This is reminiscent to the method used for front 
propagation (\cite{ES,BES}).    
We recall that the appearance of Dirac concentrations in a different area of biology 
(trait selection in evolution theory) relies also on the phase function and the viscosity 
solutions to Hamilton-Jacobi equations, \cite {DJMP, barlesP}.

In Section \ref{sec:asan} we obtain new and more precise versions of the results
of \cite{chk} by analyzing the asymptotics/rates as $\sigma \to 0$.
In Section \ref{sec:large} we present new results for large transition coefficients,
while in Section \ref{sec:vanish} we consider coefficients that may
vanish.

\section{Bounded non-vanishing transition coefficients}
\label{sec:asan}

We begin with the assumptions on the transition rates and
potentials. As far as the former are concerned we assume that
\beq \text{ there exists }\quad k>0 \quad \text{such that} \quad
\nu_{ij} \geq k>0 \qquad \text{ for all} \qquad i\neq j.
\label{as:nu1} \eeq

As far as the potentials are concerned we assume that, for all
$i=1,2,...,I,$

\beq \psi_i \in C^{2,1}(0,1), \label{as:psi1} \eeq
\beq \text{ there exists a finite collection of intervals} \quad
(J_k)_{1\leq k\leq M}\quad \text{ such that } \; \dis \min_{1\leq
i\leq I} \psi'_i
>0 \quad \text{in}\quad \bigcup J_k, \label{as:psi2} \eeq
and
\beq \max_{1\leq i\leq I}  \psi_i' >0 \qquad \text{ in} \qquad
[0,1]. \label{as:psi3} \eeq

Notice that these assumptions are satisfied by periodic potentials
with period $1/M$.

\begin{figure}[ht!]
\begin{center}
\end{center}
\vspace{-.7cm} \caption{ \label{fig:biomotor1} {\footnotesize\sc
Motor effect exhibited by the parabolic
system \fer{eq:biomotor} with two asymmetric potentials. Left: the potentials $\psi_1$, $\psi_2$.  Right: the phase functions $R_1^\sg = -\sg \ln (n_1^\sg) $, $R_2^\sg = - \sg \ln (n_2^\sg)$.  As announced in Theorem \ref{th:main},  we have $R_1^\sg \approx  R_2^\sg$ and are nondecreasing. This means that the densities are concentrated as Dirac masses at $x=0$. Here we have used $\sg=10^{-4}$.  See Figure \ref{fig:biomotor2} for another behavior. }
        }
\end{figure}

Our first result is a new and more precise version of the result in
\cite{chk}. It yields that the system \fer{eq:biomotor} exhibits a
motor effect for $\sg$ small enough and molecules are necessarily
located at $x=0$. This effect is explained by a precise asymptotic
result in the limit $\sg \to 0$.

To emphasize the dependence on the diffusion
$\sigma$, in what follows we denote, for all $i=1,...,I$, by $n^{\sg}_i$ the
solution of \fer{eq:biomotor}. Moreover, instead of \fer{eq:norm},
we use the normalization

\beq \dis \sum_{1\leq i \leq I}  n^\sg_i(0)=1.
\label{eq:norm1}
\eeq

We have:

\begin{theorem} Assume that \fer{as:nu}, \fer{as:nu1}, \fer{as:psi1},
\fer{as:psi2}, \fer{as:psi3} and \fer{eq:norm1} hold. Then, for all
$i=1,...,I$,
$$
R^{\sg}_i=- \sg \ln n^{\sg}_i{\;}_{\overrightarrow{\; \sg
\rightarrow 0 \; }}\; R \quad \text{ in } \quad C(0,1), \quad R(0)=0
\quad \text{ and }\quad R'=\dis \min_{1\leq i\leq I}  ( \psi'_i )_+ \; .
$$
\label{th:main}
\end{theorem}

In physical terms, $R$ can be seen as an effective potential for the system. To 
state the next result, we recall that throughout the paper we
denote by $\delta_0$ the usual $\delta$-function at the origin.

We have:

\begin{corollary} Assume, in addition to \fer{as:nu}, \fer{as:nu1}, \fer{as:psi1},
\fer{as:psi2} and \fer{as:psi3}, that $\dis \min_{1\leq i \leq I}
\psi'_i(0)>0$, and normalize $n^\sg_i$ by \fer{eq:norm} instead of
\fer{eq:norm1}. There exist $(\rho_i)_{1\leq i\leq I}$ such that
$$
n^\sg_i    {\;}_{\overrightarrow{\; \sg \rightarrow 0 \; }}\; \rho_i
\delta_0, \quad 
\rho_i > 0, \quad \text{and} \quad \dis \sum_{1\leq i\leq I}
\rho_i=1.
$$
\label{cor:main}
\end{corollary}

There are several possible extensions of Theorem \ref{th:main}. Here we
state one which, to the best of our knowledge, is not covered by
any of the existing results.

To formulate it, we need to introduce the following assumption on
the potentials $(\psi_i)_{1\leq i \leq I}$ which replaces
\fer{as:psi3} and allows to consider more general settings. It is:

\beq \left\{
\begin{array}{l}
\text { the set} \quad \{ x \in [0,1]: \; \dis \max_{1\leq i\leq I}
\psi_i' (x)< 0\} \quad \text { is a union of finitely  many
intervals } \quad (K_l)_{1\leq l \leq
M'},\\ \\
\text{ and } \quad (\cup J_k)^c\cap (\cup K_l)^c \quad \text { is
either a finite union of intervals or isolated points.}
\end{array}\right.
\label{as:psi4}
\eeq

We have:

\begin{theorem} Assume \fer{as:nu}, \fer{as:nu1}, \fer{as:psi1},
\fer{as:psi2}, \fer{as:psi4} and \fer{eq:norm1}. Then
$$
R^{\sg}_i=- \sg \ln n^{\sg}_i {\;}_{\overrightarrow{\; \sg \rightarrow 0
\; }}\; R \quad \text{ in } \quad C(0,1), \quad R(0)=0 \quad \text{
and } \quad
$$
$$
R' =\left\{ \begin{array}{ll}
 \dis \min_{1\leq i\leq I}  ( \psi'_i )_+  \quad \text {in} \quad \cup J_k,
\\ \\
 \dis \max_{1\leq i\leq I}  \psi'_i \quad  \text {in} \quad \cup K_l ,
\\ \\
0 \quad \text{in} \quad{\rm Int} \big((\cup J_k)^c \cap (\cup K_l)^c\big) .
\end{array} \right.
$$
\label{th:main2}
\end{theorem}

As a consequence we have:

\begin{corollary} In addition to  \fer{as:nu}, \fer{as:nu1}, \fer{as:psi1},
\fer{as:psi2}, \fer{as:psi4} and \fer{eq:norm}, assume that we have
the same number of intervals $J_k$ and $K_l$ in \fer{as:psi2} and
\fer{as:psi4} respectively, that
$0$ is the left endpoint of $J_1$  and, finally, that, for all $k=1,...,M$,
$$
\big| \int_{K_k} \; \max_{1\leq i\leq I} \psi_i'(y) dy \big| <
\int_{J_k} \; \dis \min_{1\leq i\leq I}  \psi_i'(y) dy.
$$
Then, for all $i=1,...,I$, there exist $(\rho_i)_{1\leq i\leq I}$
such that
$$
n^\sg_i  {\;}_{\overrightarrow{\; \sg \rightarrow 0 \; }}\;   \rho_i
\delta_0, \quad 
\rho_i > 0, \quad \text{and} \quad \dis \sum_{1\leq i\leq I}
\rho_i=1.
$$
\label{cor:main2}
\end{corollary}

\medskip

Other possible extensions concern coefficients that may vanish
somewhere and/or be unbounded. The former case is studied in Section
\ref{sec:vanish}. As far as the $\nu_{ij}$ being unbounded, it will be clear from
the proof of Theorem \ref{th:main}, that the coefficients  can
depend on $\sg$ as long as, for $\sg \to 0$ and all $i,j=1,...,I$,
there exists $\alpha
>0$ such that
$$
\sg \nu_{ij}^\sg \to 0 \qquad \text { and } \qquad \sigma^{-\alpha}
\nu_{ij} \to \infty .
$$
Going further in this direction leads to a different limits for
$-\ln {n_i}^\sigma$ that we study in the next Section.
\\

We continue next with the proof of Theorem \ref{th:main}. The modifications needed
to prove Theorem \ref{th:main2} are indicated at the end of this section where we also
discuss the proofs of the Corollaries.
\\

\noindent {\bf Proof of Theorem \ref{th:main}} A direct computation
shows that the $R_i^\sg$'s satisfy, for $\wt \nu_{ii}=\nu_{ii} -
\psi''_i$, the system

 \beq
\label{eq:rsigma} \left\{
\begin{array}{l}
- \sg \f{\p^2 R_i^\sg}{\p x^2} + \f{\p R_i^\sg}{\p x}^2 - \psi'_i(x)
\f{\p R_i^\sg}{\p x} + \sg  \dis \sum_{j=1}^I
\nu_{ij}e^{(R^\sg_i-R^\sg_j)/\sg} =  \sg \wt \nu_{ii} \quad \text{
in} \quad (0,1),
\\  \\
\f{\p R_i^\sg}{\p x} =\psi'_i \quad \text{ in } \{0,1\}.
\end{array}
\right. \eeq

Adding the equations of  \fer{eq:biomotor} and using \fer{as:nu}
yield the conservation law
$$
- \sg  \f{\p^2}{\p x^2} [\dis \sum_{1\leq i \leq I} n_i^\sg ] - \f{\p}{\p x}
[\dis \sum_{1\leq i \leq I} \psi'_i n_i^\sg] =0,
$$
which together with the boundary condition gives

\beq \label{eq:zeroflux} - \sg \f{\p}{\p x} \dis \sum_{1 \leq i \leq
I} n_i^\sg - \sum_{1\leq i \leq I} \psi'_i n_i^\sg =0. \eeq

Setting
$$
 \sum_{1\leq i\leq I}  n_i^\sg = e^{-S^\sg /\sg},
$$
we have
$$
\f{\p S^\sg}{\p x}= \f{\sum_i \psi'_i n_i^\sg}{\sum _i n_i^\sg},
$$
and, as a consequence, the {\it total
flux estimate}

\beq \min_{1\leq i\leq I} \psi'_i \leq \f{\p S^\sg}{\p x} \leq
\max_{1\leq i\leq I}  \psi'_i . \label{eq:sbounds} \eeq

The normalization \fer{eq:norm1} of  the $n_i^\sg$'s implies that
$S^\sg(0)=0$. As a result, there exists a $S\in C^{0,1}(0,1)$ such
that, after extracting a subsequence,
 \beq \left\{
\begin{array}{l} S^\sg {\;}_{\overrightarrow{\; \sg \rightarrow 0 \;
}}\;  S, \qquad S(0)=0, \quad \text { and}
\\ \\
\dis \min_{1\leq i\leq I}  \psi'_i  \leq  \f{\p S}{\p x} \leq
\max_{1\leq i\leq I}  \psi'_i \quad \text{in } [0,1] .
\end{array} \right.
\label{eq:ests}
\eeq

Next we obtain bounds on the $R^\sigma_i$'s, which are independent
of $\sigma$, and imply their convergence as $\sigma \to 0$. This is
the topic of the next Lemma which we prove after the end of the
ongoing proof.

\begin{lemma} For each $i=1,...,I$ there exists a positive  constant
$C_i= C_i(\psi'_i, \sg \nu_{ii}, \sigma \psi''_i))$ such that
$$
|R_i^\sg | + |\f{\p R_i^\sg}{\p x}| \leq C_i \quad \text{ in}\quad
[0,1],
$$
Moreover, for all $i=1,...,I$,
$$
R_i^\sg {\;}_{\overrightarrow{\; \sg \rightarrow 0 \; }}\;   R=S ,
\qquad \text { in } C([0,1]).
$$
\end{lemma}

We obtain next the Hamilton-Jacobi satisfied by the limit $R=S$. The
claim is that the limit is a viscosity solution (see, for instance,
\cite{barles, CIL}) of

\beq \label{eq:hj} \left| \f{\p R}{\p x}\right|^2 + \dis\max_{1\leq
i\leq I} \; [- \psi'_i \f{\p R}{\p x}]=0 \qquad \text{ in} \quad
(0,1). \eeq 

We do not state the boundary condition because we do not
use them. It can, however, be proved that $R$ satisfies
$$
 \f{\p
R}{\p x}\leq \dis\max_{1\leq i\leq I} \psi'_i \quad \text{at}\quad
x=0 \quad \text{ and} \quad  \f{\p R}{\p x}\geq \dis\min_{1\leq
i\leq I} \psi'_i \quad \text{ at} \quad  x=1. $$

We begin with the subsolution property. Letting $\sigma \to 0$  in
the inequality
$$
- \sg \f{\p^2 R_i^\sg}{\p x^2} + \left| \f{\p R_i^\sg}{\p x}
\right|^2- \psi'_i \f{\p R_i^\sg}{\p x} \leq  \sg \wt \nu_{ii},
$$
gives, for all $i=1,...,I$,
$$
\left| \f{\p R}{\p x}\right|^2 - \psi'_i \f{\p R}{\p x} \leq 0.
$$
\\

To prove that $R$ is a supersolution of \fer{eq:hj} we observe that
function $R^\sg = \dis \min_{1\leq i\leq I} R_i^\sg$ satisfies the
inequality
$$ - \sg \f{\p^2 R^\sg}{\p x^2} + \left|\f{\p R^\sg}{\p x}\right|^2 +
\max_{1\leq i\leq I} \; [- \psi'_i(x) \f{\p R^\sg}{\p x} ]+ \sg
\sum_{i,j=1}^I \nu_{ij} \geq  \sg \min_i (\wt \nu_{ii}).
$$

Letting again $\sigma \to 0$, we find that $R=S= \dis \lim_{\sg \to 0}
R^\sg$ satisfies
$$
\left | \f{\p R}{\p x} \right|^2 + \max_{1\leq i\leq I} \; [ -
\psi'_i(x) \f{\p R}{\p x} ]\geq 0.
$$
\\

We obtain now the formula for $R$. To this end, observe first that, since
$ \dis
\lim_{\sg \to 0} R^\sg=R=S$,
letting $\sigma \to 0$ in \fer{eq:sbounds} yields

 \beq \label{eq:lim1} \min_{1\leq i\leq
I} \psi'_i \leq \f{\p R}{\p x} \leq \max_{1\leq i\leq I} \psi'_i .
\eeq

Next we show that, in the viscosity sense, \beq \label{eq:p} \f{\p
R}{\p x} \geq 0. \eeq

Indeed for a  test function $\Phi$, let $x_0 \in (0,1)$ be the
maximum of $R-\Phi$, i.e., $(R-\Phi)(x_0) = \dis \max_{0\leq x\leq
1}(R - \Phi)(x)$ and assume that
$$\Phi'(x_0)<0.$$

Applying the viscosity subsolution criterion to \fer{eq:lim1}, then
implies that
$$\Phi'(x_0)- \max_i \psi'_i(x_0) \geq 0.$$

This, however, contradicts the inequality
$$\dis \max_{1\leq i\leq I}
\psi'_i(x_0)  >0$$ that follows from the assumption \fer{as:psi2}.

Combining \fer{eq:lim1} and \fer{eq:p} we get
 \beq \label{eq:lim2} \min_{1\leq
i\leq I} (\psi'_i )_+ \leq \f{\p R}{\p x} \leq \dis\max_{1\leq i\leq
I} \psi'_i . \eeq

Finally, given a test function $\Phi$, let $x_0 \in (0,1)$ be such
that $(R- \Phi)(x_0) = \dis \max_{0\leq x\leq 1} (R- \Phi)$ and
assume that $$\Phi'(x_0)>0.$$ 

Again by the viscosity criterion we
must have
$$\Phi'(x_0)- \dis \min_{1\leq i\leq I} \psi'_i(x_0) \leq 0,$$
and, hence, in the viscosity sense,
\beq \label{eq:lim3}
 \f{\p R}{\p x} \leq (\min_{1\leq i\leq I} \psi'_i )_+  \qquad 
 \text{if }\;  \f{\p R}{\p x} >0.
\eeq

This concludes  the proof of the formula in the claim.
 \qed

We return now to the

\noindent {\bf Proof of Lemma 3.1} For the Lipschitz estimate,
observe that, at any extremum point $x_0$ of 
$\f{\p R_i^\sg}{\p x}$, we have
$\f{\p^2 R_i^\sg}{\p x^2} =0$. Evaluating the equation at $x_0$, we get
$$
\left| \f{\p R_i^\sg}{\p x}\right|^2 \leq  \psi'_i \f{\p
R_i^\sigma}{\p x} + \sg \wt \nu_{ii}.
$$

As a consequence, at $x_0$ we have
$$
\left| \f{\p R_i^\sg}{\p x} \right| \leq \max_{0\leq x \leq 1}
\psi'_i +\sqrt{ \sg \wt \nu_{ii}}.
$$

To identify the limit of $ \min_{1\leq j \leq I}  R_j^\sg$ notice
that  the inequality $$ n_i^\sg \leq \dis \sum_{1\leq j \leq I}
n_j^\sg \leq I \max_j n_j^\sg$$ gives
$$
- \sg \ln I + \min_{1\leq j \leq I}  R_j^\sg  \leq S^\sg  \leq
R_i^\sg,$$ and  thus
$$S^\sg  \leq  \min_{1\leq i\leq I} R_i^\sg.
$$

Consequently, we have  the uniform convergence
$$
\min_{1\leq i\leq I} R_i^\sg  {\;}_{\overrightarrow{\; \sg
\rightarrow 0 \; }}\;  S.
$$

To prove the claim about the limit of the $R^{\sigma}_i$ we observe
that summing over $i$ the equations of \fer{eq:rsigma} yields
$$
\sg \sum_{i,j=1}^I \nu_{ij} ( \f{(R_j^\sg-R_i^\sg)_+}{\sg})^2 \leq 2
\sg \sum_{i,j=1}^I \nu_{ij}e^{(R^\sg_i-R^\sg_j)/\sg} \leq 2(\sg
\sum_{1\leq i \leq I} \wt  \nu_{ii}+
 2 \sg \f{\p^2 \sum_i R_i^\sg}{\p x^2}+ \sum_{1\leq i \leq I}  \psi'_i \f{\p R_i^\sg}{\p x}) .
$$

Integrating in $x$ and using the gradient estimates, we find that
$$
\dis \sum_{i,j=1}^I \;  \int_0^1 ( R_j^\sg-R_i^\sg)^2 =\f 1 2 \dis
\sum_{i,j=1}^I \;  \int_0^1 ( R_j^\sg-R_i^\sg)_+^2 \leq C \sg.
$$

Together with the uniform gradient estimate on $R_i^\sg$ and the
uniform bound on $\min_{1\leq j \leq I}  R_j^\sg$, we deduce that
$$
R_i^\sg  {\;}_{\overrightarrow{\; \sg \rightarrow 0 \; }}\;  R=S \in
C^{0,1}(0,1).
$$ \qed

We continue with the

\noindent {\bf Proof of Corollary \ref{cor:main}} The normalization
\fer{eq:norm} amounts to adding a constant to the $R_i$.  The
exponential behavior of $n_i^\sg$, with an increasing $R_i^\sg$
(from Theorem \ref{th:main}), yields  that the $n_i^\sg$'s converge,
as $\sigma \to 0$,  to $0$ uniformly on intervals $[\e, 1]$ with
$\e>0$. Moreover,  $R(0)=0$. The result follows with $\rho_i \geq
0$. If $\rho_i=0$ for some $i=1,...,I$, then, letting $\sigma \to 0$
in  \fer{eq:biomotor}, gives, in the sense of distributions, that
$$
0=\dis{\sum_{j\neq i}}\;  \nu_{ij} n_j.
$$

But then all the $\rho_j$ must vanish, which is impossible with the
normalization of unit mass. \qed

We present now a brief sketch of the proof of Theorem
\ref{th:main2}. Since it follows along the lines of the proof of
Theorem \ref{th:main}, here we only point out  the differences.

We have:

\noindent {\bf Proof of Theorem \ref{th:main2}} The Lipschitz
estimates,  the passage in the limit and the identification of the
limiting Hamilton-Jacobi equation in the Theorem \ref{th:main} did
not depend on the assumption \fer{as:psi3}, hence, they hold true
also on the case at hand. The final arguments of the proof of
Theorem \ref{th:main2} also identify the limit on the set $(\cup
K_l)^c$. On the set $\cup K_l$ we already know from \fer{eq:lim1}
that $R'$ is less than the claimed value, and thus it is negative.
We conclude the equality by using the Hamilton-Jacobi equation.
Indeed in this situation we know that
$$
\dis\max_{1\leq i\leq I} \; [- \psi'_i \f{\p R}{\p x}]=
 - \f{\p R}{\p x} \dis\max_{1\leq i\leq I} \; \psi'_i.
$$\qed

We conclude the section with the proof Corollary 2.4, which is
simply a variant of the one for Corollary 2.2. We have:

\noindent{\bf Proof of Corollary 2.4}  The assumption on $\cup J$
asserts that $R$ is increasing on $\cup J$. Then it may decrease
but, for $x>0$, $R(x)>R(0)$. With the unit mass normalization, this
means that $R(0)=0$ as before and the convergence result holds as
before. \qed

\section{Large transition coefficients} 
\label{sec:large}

\begin{figure}[ht!]
\begin{center}
\end{center}
\vspace{-.7cm} \caption{ \label{fig:biomotor2} {\footnotesize\sc
Motor effect exhibited by the parabolic
system \fer{eq:biomstrong} with large transition coefficients. The figure depicts the phase functions $R_1^\sg$, $R_2^\sg$.  As announced in Theorem \ref{th:strong}, we have $R_1^\sg \approx  R_2^\sg$ and can decrease slightly. Here we have used $\sg=5 \; 10^{-3}$.}
        }
\end{figure}

In this section we consider transition coefficients normalized by
$1/\sg$. For the sake of simplicity we take $I=2$. This allows for
explicit formulae. The equations for larger systems, i.e.,  $I>3$,
are more abstract. The system \fer{eq:biomotor} is replaced by
\beq
\left\{
\begin{array}{ll}
-\sg \f{\p^2}{\p x^2}  n^{\sg}_1 - \f{\p}{\p x}( \nabla \psi_1\;
n^{\sg}_1) +\f 1 \sg \nu_{1} n^{\sg}_1  = \f 1 \sg \nu_{2} n^{\sg}_2
\quad \text{ in} \quad (0,1),
\\ \\
-\sg \f{\p^2}{\p x^2}  n^{\sg}_2 - \f{\p}{\p x}( \nabla \psi_2\;
n^{\sg}_2) +\f 1 \sg \nu_2 n^{\sg}_2  =  \f 1 \sg \nu_{1} n^{\sg}_1
\quad \text{ in} \quad (0,1),
\\ \\
\sg \f{\p}{\p x}  n^{\sg}_i + \nabla \psi_i \; n^{\sg}_i=0 \quad
\text{ in } \quad  \{0,1\} \quad \text{ for }  \quad i=1, 2.
\end{array}
\right.
\label{eq:biomstrong}
\eeq

As before we assume that

\beq \label{as:pos1} n^{\sg}_i >0 \quad \text { in } [0,1] \quad
\text {for} \quad i=1, 2. \eeq

The result is:

\begin{theorem} Assume \fer{as:nu},  \fer{as:nu1},  \fer{as:psi1}, \fer{as:psi2} and
\fer{as:psi3} and consider the solution $(n^\sg_1, n^\sg_2)$ to
\fer{eq:biomstrong}  normalized  by $ n^\sg_1(0)+n^\sg_2(0)=1$.
Then, as $\sg \to 0$ and $i=1,2$,
$$
R^\sg_i=- \sg \ln n^\sg_i{\;}_{\overrightarrow{\; \sg \rightarrow 0
\; }}\;  R \quad \text{ in } C(0,1), \qquad R(0)=0, \quad\text{and}
$$
$$
 R' \geq  \left\{ \begin{array}{ll}
           \dis \min_{1\leq i\leq I}  \psi'_i &  \text{ on } \; \cup J_l ,
           \\ \\
           -\sqrt{k} &  \text{ on } \; (\cup J_l)^c.
 \end{array}
\right.
$$
\label{th:strong}
\end{theorem}

The corollary below follows from Theorem \ref{th:strong} in a way
similar to the analogous corollaries in the previous section. Hence,
we leave the details to the reader.

\begin{corollary} In addition to the assumptions of Theorem \ref{th:strong}, suppose
that $0 \in J_1$, the potentials  are small enough so that
$$
\sqrt{k} \; | K |  < \int_{\sup J} \; \dis \min_{1\leq i\leq 2}
\psi_i'(y) dy,
$$
and $(n^\sg_1,n^\sg_2)$ is normalized by \fer{eq:norm}. There exist
$\rho_1,\rho_2 >0$ such that $\rho_1 + \rho_2 =1$ and, as $\sigma
\to 0$ and for $i=1,2$,
$$
n^\sg_i  {\;}_{\overrightarrow{\; \sg \rightarrow 0 \; }}\;   \rho_i
\delta_0.
$$
\label{cor:strong}
\end{corollary}

We present next a sketch of the proof of Theorem \ref{th:strong} as
most of the details follow as in the previous theorems.

\noindent {\bf Proof of Theorem \ref{th:strong}}
 The total flux and
Lipschitz estimates follow as before. The main new point is the
limiting Hamilton-Jacobi equation which is more complex. We
formulate this as a separate lemma below. Its proof is based on the
use of perturbed test functions. We refer to  \cite{BES} for the
rigorous argument in a more general setting.

\begin{lemma} The uniform in $[0,1]$ limit $R$, as $\sigma \to 0$,  of the $R_i^\sg$ satisfies
the Hamilton-Jacobi equation \beq \label{eq:hjstrong} H\big( \f{\p
R}{\p x}, x \big) =0, \quad \text{ in } \quad (0,1) \eeq with \beq
\label{eq:ham} H(p,x)= \f 1 2 [ \beta_1+ \beta_2 +\sqrt{( \beta_1+
\beta_2)^2 -4 ( \beta_1 \beta_2 - \nu_1 \nu_2 }] , \eeq where, for
$i=1,2,$
$$
\beta_i= p^2 - \psi'_i \; p -\nu_{i}.$$
\end{lemma}

The formula for $R'$ follows from the above  Lemma by analyzing the
solutions to the Hamilton-Jacobi equation as before. On the set
$\cup J_l$ the answer follows from the bounds \fer{eq:sbounds}. On
the set $(\cup J)^c$ the argument is more elaborate. Using that $R'$
is a subsolution, we get
$$
\beta_1 \beta_2 - \nu_1 \nu_2 \geq 0 \quad \text{ and} \quad \beta_1
+ \beta_2 \leq 0.
$$

Therefore both $\beta_1$ and $\beta_2$ are nonpositive and thus
$$
( R')^2 - \psi'_i R' -\nu_i \leq 0.
$$

On the other hand we know that on $(\cup J_l)^c$ one of the
potentials -- for definiteness say $\psi _1$ --  satisfies $\psi'_1
>0$,
hence, always in $(\cup J_l)^c$,
$$
R' \geq \f 1 2 [\psi'_1(x) -\sqrt{(\psi'_1)^2+4 \nu_1}  ] \geq
\sqrt{ \nu_1}.
$$

The inequalities for $R'$ are now proved. \qed

\section{Vanishing transition coefficients}
\label{sec:vanish}

We focus here to the case where the transition
coefficients $(\nu_{ij})_{1 \leq i,j \leq I}$ may  vanish at either some
points or, in fact, on large sets. In this situation, we assume
that

\beq \left\{\begin{array}{l} \text {for each $j=1,...,I$,} \quad
\psi'_j <0 \quad \text { on a finite collection of intervals} \quad
(K_j^\al)_{1\leq \al \leq A_j} \quad \text {and}
\\ \\
\text { for all $j=1,...,I$ and  $\al = 1,...,A_j$, there exists
 $i \in \{= 1,...,I\}$  such that}\\ \\
\text { $\psi'_i \geq 0$ on $K_j^\al$}, \quad \text { and, in a left
neighborhood  of the right endpoint of $K_j^\al$}, \quad \nu_{ij}>0
.
\end{array} \right.
\label{as:nu3}
\eeq
To go for weaker assumptions would face the completely decoupled case (when $\nu$ vanishes) and the motor effect does not occur. 
\\

We have:
\begin{theorem} Assume \fer{as:psi1},\fer{as:psi2},\fer{as:psi3}, \fer{as:nu3}
and normalize the solution $(n^\sg_i)_{1\leq i \leq I}$ to
\fer{eq:biomstrong} by \fer{eq:norm1}.
For $i=1,...,I$, let $R_i^\sigma =-\ln n_i^\sigma$.
Then, as $\sg \to 0$,
$$
\text{either }  \quad R_i^\sg {\;}_{\overrightarrow{\; \sg
\rightarrow 0 \; }}\;  R_i \quad \text{ in } C(0,1),
\quad \text { or } \quad  R_i^\sg {\;}_{\overrightarrow{\; \sg
\rightarrow 0 \; }}\;  \infty  \quad \text{ uniformly in} \quad
[0,1].
$$
Moreover, the function $R= \dis \min_{1\leq i \leq I}  R_i$
satisfies
$$
R(0)=0, \quad R'\geq 0 \quad \text{ and } \quad R' = \dis
\min_{1\leq i \leq I} \psi'_i \quad \text{ on } \quad \cup J_l.
$$
\label{th:vanish}
\end{theorem}

We also have:

\begin{corollary} In addition to the assumptions of Theorem \ref{th:vanish}, suppose that
$0\in J_1$. For $i=1,,,I$, there exist $\rho_i \geq 0$ such that
$\sum_{1\leq i \leq I}  \rho_i =1$, and, as $\sigma \to 0$,
$$
n^\sg_i {\;}_{\overrightarrow{\; \sg \rightarrow 0 \; }}\;  \rho_i
\delta_0.
$$
\end{corollary}

The direct conclusion of Theorem \ref{th:vanish} is simply that
$$
\sum_{1\leq i \leq I} n^\sg_i  {\;}_{\overrightarrow{\; \sg
\rightarrow 0 \; }}\; \delta_0 $$.

The corollary follows from the fact that, for all $i=1,...,I,$
$n_i^\sg \geq 0$. We do not know whether in this context each
$\rho_i$ is positive. To get this, we need to assume something more
like, for example, $\nu_{ij} (0) >0$ for all $i,j=1,...,I$.

We conclude with a brief sketch of the

\noindent {\bf Proof of Theorem \ref{th:vanish}}. The total flux and
Lipschitz estimates follow as before. A careful look at the proof of
the convergence part of Theorem \ref{th:main} shows that either the
$R_i^\sg$'s  blow up or they are uniformly bounded and, hence,
converge uniformly in $(0,1)$ to a subsolution of
$$
|R_i'|^2 -\psi'_i R'_i \leq 0.
$$

It then follows that
$$
R'_i \geq 0  \text{ on } \left( \dis \cup_\al K_i^\al \right)^c
\quad \text {and } R'_i \leq \psi'_i(x)   \text{ on } \dis \cup_\al
K_i^\al .
$$

The final step is to prove that
$$
R(x)=\dis \min_{i\in L(x)} R_i(x) \quad \text { in }  \quad L(x)=\{
i, \; \psi'_i(x) \geq 0\}.
$$

This follows as before. We leave the details to the reader.\qed

%
%

\bigskip
\bigskip

\noindent ($^{1}$) Ecole Normale Sup\'erieure \\
             DMA, UMR8553 \\
             45 rue d'Ulm, 75230 Paris \\
             France \\
             and Institut Universitaire de France
             \\
email: benoit.perthame@ens.fr \\ \\
($^{2}$)  Department of Mathematics \\
             The University of Texas at Austin \\
             Austin, TX 78712 \\
             USA \\
email: souganid@math.utexas.edu\\ \\
($^{3}$)  Partially supported by the National Science Foundation.

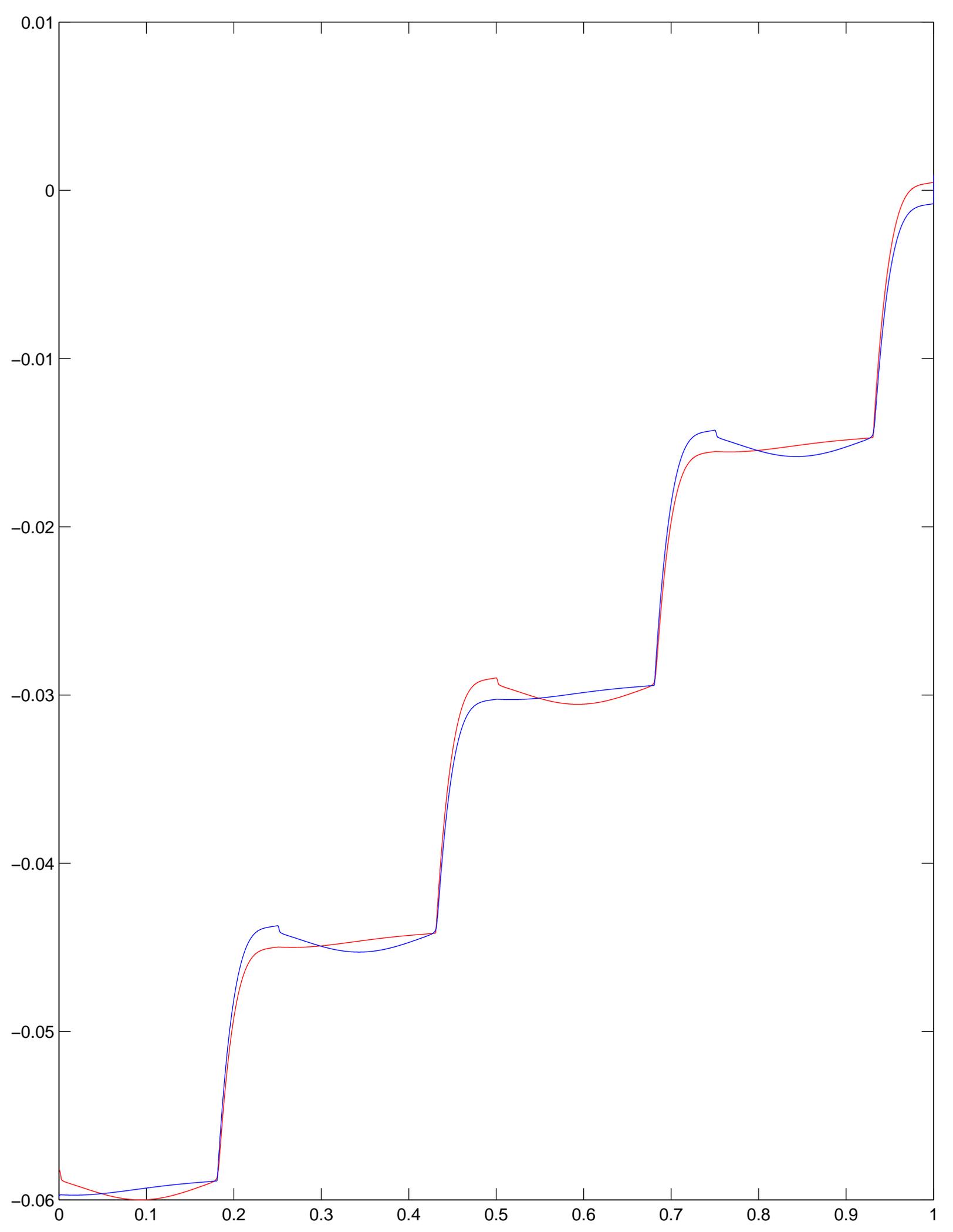
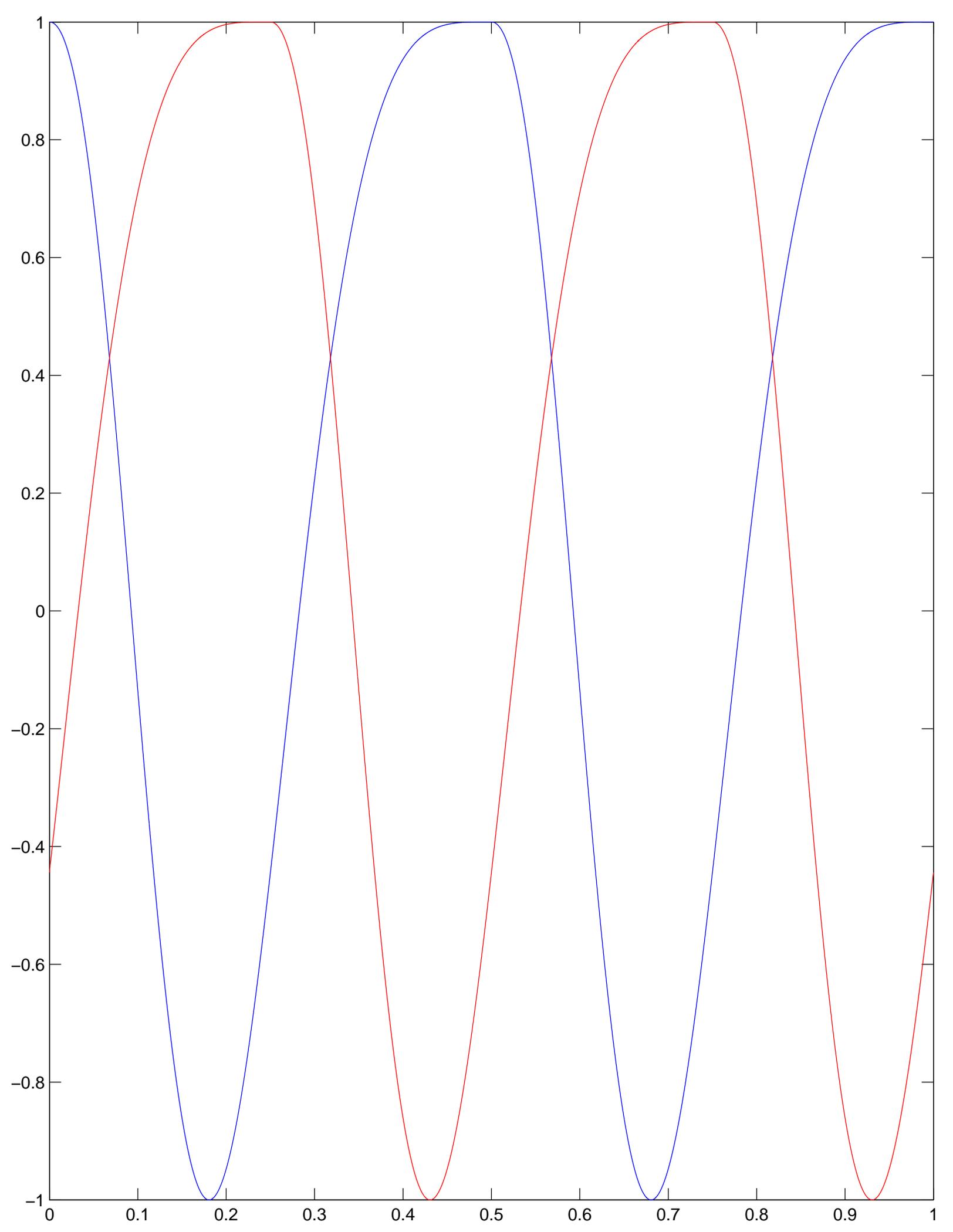
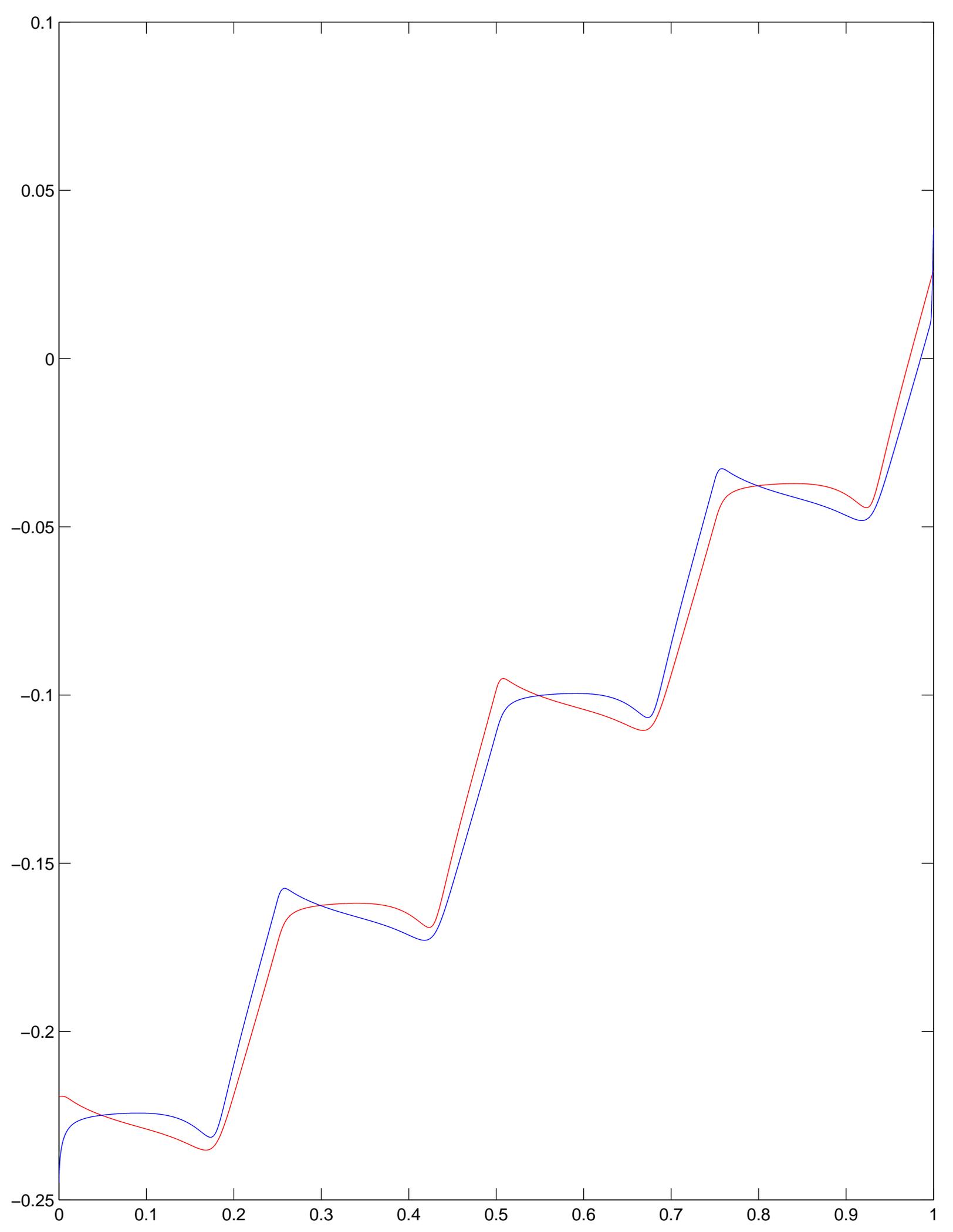
\end{document}